\documentstyle[12pt]{article}
\input{amssym.def}
\input{amssym}

\def\be{\begin{equation}}
\def\ee{\end{equation}}
\def\ba{\begin{eqnarray}}
\def\ea{\end{eqnarray}}
\def\lb{\label}
\def\nin{\noindent}
\def\tr{{\rm Tr}}

\begin{document}
\begin{center}
\vskip .5cm

{\LARGE{\bf On $R$-matrix representations of Birman-Murakami-Wenzl algebras}}

\vskip .6cm

{\large {\bf A. P. Isaev}}

\vskip 0.2 cm

Bogoliubov Laboratory of Theoretical Physics, JINR\\
141980 Dubna, Moscow Region, Russia

\vskip .4cm

{\large {\bf O. V. Ogievetsky\footnote{On leave of absence from P. N.
Lebedev Physical Institute, Leninsky Pr. 53, 117924 Moscow, Russia}}}

\vskip 0.2 cm

Center of Theoretical Physics, Luminy \\
13288 Marseille, France

\vskip .4cm

{\large {\bf P. N. Pyatov}}

\vskip 0.2 cm

Bogoliubov Laboratory of Theoretical Physics, JINR\\
141980 Dubna, Moscow Region, Russia

\end{center}

\vspace{.5 cm}
\nin {\sc Abstract.} {\sf We show that to every local representation of
the Birman-Murakami-Wenzl algebra defined by a skew-invertible $R$-matrix 
$\hat R\in \mbox{Aut}(V^{\otimes 2})$ one can associate pairings 
$V\otimes V\rightarrow {\Bbb C}$ and $V^*\otimes V^*\rightarrow {\Bbb C}$, 
where $V$ is the representation space. Further, we investigate conditions 
under which the corresponding quantum group is of $SO$ or $Sp$ type.}

\vskip .5cm
\hfill {\sl To the memory of Andrei Nikolayevich Tyurin} 
\vskip .6cm

\setcounter{equation}0

Let $G$ be either orthogonal or symplectic Lie group, ${\goth g}$ its Lie 
algebra and ${U}_q({\goth g})$ the corresponding quantum group (i.e., the 
quantized universal enveloping algebra \cite{D,FRT}). Denote by $V$ the 
space of the vector representation of $G$ or ${U}_q({\goth g})$.

In \cite{Br} R. Brauer constructed centralizers $\mbox{End}_{G}(V^{\otimes n})$
of the action of $G$ on tensor powers of the vector representation.
He introduced a one-parametric family of algebras ${Br}_n(x)$; for certain 
values of the parameter $x=x_{G}$~\footnote{Values $x_{G}$ and $\nu_{G,q}$ 
(introduced in a paragraph below) depend essentially on a particular choice 
of the group $G$.}, the algebras ${Br}_n(x_{G})$ possess representations 
${  B  r}_n(x_{  G}) \rightarrow \mbox{End}(V^{\otimes n})$ 
commuting with the action of $G$. These representations are generated by 
the permutation $P\in \mbox{Aut}(V^{\otimes 2})$: $P(u\otimes v) = 
v\otimes u , \;\forall\; u,v\in V$ and the operation related to the 
$G$-invariant pairing $g:\ V\otimes V\rightarrow {\Bbb C}$.

In case of ${U}_q({\goth g})$ the role of the Brauer centralizer algebras 
is played by a two parametric family of algebras ${W}_n(q,\nu)$ introduced 
independently by J. Murakami \cite{M} and by J. Birman and H. Wenzl \cite{BW}.
Centralizers $\mbox{End}_{{U}_q({g})}(V^{\otimes n})$ are then realized by 
specific representations of the Birman-Mirakami-Wenzl algebras  
${W}_n(q,\nu_{{G},q}) \rightarrow \mbox{End}(V^{\otimes n})$
which are generated by $q$-analogs of permutations called the R-matrices.
(all necessary definitions are given below; references on the topic can be 
found in a review book \cite{ChP}).   

In the present paper we study an inverse problem. Given an R-matrix 
representation of the Birman-Murakami-Wenzl algebra we find conditions under 
which the associated quantum group can be called an orthogonal or
symplectic one. More precisely, we prove that for any R-matrix which 
generates representations of algebras ${W}_n(q,\nu)$ on spaces $V^{\otimes n}$
one can construct a unique, up to a multiplicative constant, nondegenerate 
pairing on the space $V$ (see Theorem and Proposition 2). We further describe
conditions under which this pairing is invariant (see Proposition 3 and  
comments after it). So in the $q$-case the information on both the
permutation and the pairing is advantageously encoded in a single R-matrix.
\bigskip

An algebra ${W}_{n}(q,\nu)$ is a $(2n-1)!!$ dimensional quotient of the group 
algebra of the braid group ${\Bbb C}{\cal B}_n$. It is given in terms of 
generators $\{e_i, \kappa_i\}_{i=1}^{n-1}$ and relations \cite{W}
\ba
\lb{braid}
&&
e_i e_{i+1} e_i = e_{i+1} e_i e_{i+1}\ , \qquad
e_i e_j = e_j e_i\ , \qquad \;\; |i-j|>1\ ,
\\[4mm]
\lb{bmw1}
&&
e_i^2 = 1 + \lambda ( e_i - \nu \kappa_i )\ , \qquad
\lambda := q-q^{-1}\ ,
\\[2mm]
\lb{bmw2}
&&
e_i \kappa_i = \kappa_i e_i = \nu \kappa_i\ ,
\\[2mm]
\lb{bmw4}
&&
\kappa_{i+1}  e_i \kappa_{i+1} = \nu^{-1} \, \kappa_{i + 1} \ ,
\qquad
\kappa_{i+1}  e^{-1}_i \kappa_{i+1} = \nu \, \kappa_{i + 1} \ .
\ea
Here eqs. (\ref{braid}) define the Artin presentation of the braid group 
${\cal B}_n$ and  relations (\ref{bmw1})-(\ref{bmw4}) extract an
appropriate quotient. The domains of the algebra parameters 
$q\in {\Bbb C}\backslash \{0,\pm 1\}$ and
$\nu\in {\Bbb C}\backslash \{0, q, - q^{-1}\}$ are chosen in such a way that 
the elements $\kappa_i$ can be expressed in terms of $e_i$ 
\ba
\lb{def-kappa}
\kappa_i = \lambda^{-1}\nu^{-1}(q-e_i)(q^{-1}+e_i) =\lambda^{-1}\bigl(
e^{-1}_i-e_i+\lambda\bigr)\, ,
\ea
and satisfy relations
\be
\label{knorm}
\kappa_i^2 = \mu \, \kappa_i 
\ee
with a nonzero coefficient
$\mu:={\lambda^{-1}\nu^{-1}(q - \nu)(q^{-1} + \nu)}$.~\footnote{Defining 
relations for Brauer centralizer algebra ${Br}_n(x)$ follow from relations 
(\ref{braid})-(\ref{knorm}) in a limiting case 
$\nu =q^{1-x}$, ~$q\rightarrow 1$; note that in this limit the generators 
$\kappa_i$ and $e_i$ become independent.}

Note that for generic values of $\nu$ and $q$ the set of defining relations 
(\ref{braid})--(\ref{bmw4}) is not a minimal one. To show this, start with 
the first one of eqs.(\ref{bmw4}) and multiply both sides by 
$e_i^{-1}e_{i+1}^{-1}$ from the right
\be
\lb{bmw3a}
\kappa_{i + 1} \kappa_i =
\kappa_{i+1} e^{-1}_i e^{-1}_{i + 1} \ .
\ee
Then, multiplying by $\lambda\kappa_{i+1}$ from the right and performing 
straightforward transformations we get
$$
\begin{array}{c}
\lambda \nu^{-1} \kappa_{i+1} e^{-1}_{i} \kappa_{i+1} =
\lambda  \kappa_{i+1} \kappa_i  \kappa_{i+1} = \\[2mm]
= \kappa_{i+1}(e_{i }^{-1} - e_{i } + \lambda) \kappa_{i+1} =
\kappa_{i+1} ( e^{-1}_{i } \kappa_{i+1}  - \nu^{-1} + \lambda \mu )  \; ,
\end{array}
$$
wherefrom it follows that $(\lambda\nu^{-1} -1) (\kappa_{i+1} e^{-1}_{i } 
\kappa_{i+1} - \nu \kappa_{i+1}) = 0$. Thus, in case $\nu\neq\lambda$ the 
two relations in (\ref{bmw4}) are algebraically dependent.

In the sequel we shall use a few more relations for the generators $e_i$ 
and $\kappa_i$:
\ba
\lb{bmw3}
&&
\kappa_{i} \kappa_{i \pm 1} =  \kappa_{i} e_{i\pm 1} e_{i}\ ,
\\ [2mm]
\lb{bmw5}
&&
\kappa_i \kappa_{i \pm 1} \kappa_i = \kappa_i\ ,
\\[2mm]
\lb{bmw4b}
&&
\kappa_i e_{i+1}^{\pm 1} \kappa_i = \nu^{\mp 1} \kappa_i\ .
\ea
All these equalities follow easily from the defining relations 
(\ref{braid})--(\ref{bmw4}).
\bigskip

We are aiming to study a specific family of representations of algebras 
${W}_n(q,\nu)$, the so-called local (or R-matrix) representations.
We are using a compact matrix notation of \cite{FRT}.
Necessary explanations are given below.

Let $V$ be a finite dimensional vector space. We label the component 
spaces $V$ of the tensor power $V^{\otimes n} = V\otimes V\otimes \dots
\otimes V$ from left to right in the ascending order starting from 1.
For any element $X\in \mbox{End}(V^{\otimes 2})$ and for all
$1\leq k\neq l\leq n$ the symbol  $X_{kl}$ stands for an element of 
$\mbox{End}(V^{\otimes n})$ whose action differs from the identity only
on the tensor product of the $k$-th and $l$-th component spaces where it 
coincides with $X$. In case $l=k+1$ a concise notation $X_k := X_{k\, k+1}$ 
is often applied. The symbols $I$ and $P$ are reserved for the identity and 
the permutation operators respectively.

An element  ${\hat R}\in {\rm Aut}(V^{\otimes 2})$ is called an R-matrix if 
it satisfies the so called Yang-Baxter equation 
\be
\lb{YBE}
{\hat R}_1\, {\hat R}_2\, {\hat R}_1 \,=\,
{\hat R}_2\, {\hat R}_1\, {\hat R}_2\, .
\ee

With any R-matrix $\hat R$ one associates a family of representations of
the braid groups ${\cal B}_n$,
$\rho^R_n: {\cal B}_n \rightarrow {\rm Aut}(V^{\otimes n})$,
$n=1,2,\dots,$
defined on the generators by $\rho^R_n(e_i) := {\hat R}_i$, $i=1,2,\dots ,n-1$.

An R-matrix $\hat R$ whose minimal polynomial is cubic and which induces 
representations $\rho_n^R$ of the quotient algebras ${W}_n(q,\nu)$ 
(for some values of $q$ and $\nu$) is called an R-matrix of BMW type. 

An operator ${\hat R}\in \mbox{Aut}(V^{\otimes 2})$ (not necessarily an 
R-matrix) is called skew invertible iff there exists some 
$\hat{\Psi}\in {\rm End}(V^{\otimes 2})$, called the skew inverse of $\hat R$,
such that relations
\be
\lb{skew}
\tr_{(2)} (\hat{R}_{12}\hat{\Psi}_{23}) =
\tr_{(2)}(\hat{\Psi}_{12}\hat{R}_{23}) = P_{13}\, 
\ee
are satisfied. Here  the subscript $i$ in the notation of trace $\tr_{(i)}$
indicates the label of the space where the trace is evaluated. 

Denote
$$
C_2:= \tr_{(1)}\hat{\Psi}_{(12)}\ , \qquad D_1:= \tr_{(2)}\hat{\Psi}_{12} \ .
$$
As a direct consequence of the definitions  one has
\ba
\lb{trC-R}
{\tr}_{(1)}C_1 \hat{R}_{12}\  =\ I_2\ ,\quad
{\tr}_{(2)}D_2 \hat{R}_{12}\ =\ I_1\ .
\ea
In what follows while referring to relations 
(\ref{braid})--(\ref{bmw4b}) we always imply their images 
in the R-matrix representations, that is 
\be
\lb{imply}
e_i\mapsto {\hat R}_i\, , \qquad
\kappa_i\mapsto {\hat K}_i:=\lambda^{-1}\nu^{-1}(qI-\hat{R}_i)
(q^{-1}I + \hat{R}_i)\, .
\ee

\nin
{\bf Proposition~1~}\cite{R,O}.
{\it For a skew invertible R-matrix $\hat R$ the following
relations hold
\ba
\lb{psi-C}
C_1\ \hat{\Psi}_{12}\ =\ \hat{R}_{21}^{-1}\ C_2\ , &&
\hat{\Psi}_{12}\ C_1\ =\ C_2\ \hat{R}_{21}^{-1}\ ,
\\[2mm]
\lb{psi-D}
D_2\ \hat{\Psi}_{12}\ =\ \hat{R}_{21}^{-1}\ D_1\ , &&
\hat{\Psi}_{12}\ D_2\ = \ D_1\ \hat{R}_{21}^{-1}\ .
\ea
}
\vspace{-5mm}

\nin
{\bf Proof.~} First, we rewrite conditions (\ref{braid}) for the 
R-matrix $\hat R$ in the form
$
\hat{R}_{12}^{\pm 1}\ \hat{R}_{23}\ \hat{R}_{12}^{\mp 1}\ =\
\hat{R}_{23}^{\mp 1}\ \hat{R}_{12}\ \hat{R}_{23}^{\pm 1}\ .
$
Multiplying by $\hat{\Psi}_{01}\hat{\Psi}_{34}$~\footnote{
Here it is suitable to label component spaces in $V^{\otimes n}$ 
starting from 0.} and taking traces in spaces with labels 1 and 3 we get
$$
\tr_{(1)}(\hat{\Psi}_{01}\ \hat{R}_{12}^{\pm 1}\ P_{24}\ \hat{R}_{12}^{\mp 1})\ =\
\tr_{(3)}(\hat{\Psi}_{34}\ \hat{R}_{23}^{\mp 1}\ P_{02}\ \hat{R}_{23}^{\pm 1})\ .
$$
Next, evaluating trace in space 0 or 4 we get four equalities
\be
\lb{tr-04}
\tr_{(1)}(C_1\ \hat{R}_{12}^{\pm 1}\ P_{24}\ \hat{R}_{12}^{\mp 1})\ =\
C_4\ I_2\ ,
\quad
\tr_{(3)}(D_3\ \hat{R}_{23}^{\mp 1}\ P_{02}\ \hat{R}_{23}^{\pm 1})\ =\
D_0\ I_2\ ,
\ee
which can be further transformed to
\be
\lb{tr-04a}
\tr_{(1)}(C_1\ \hat{R}_{12}^{\pm 1}\  \hat{R}_{14}^{\mp 1})\ =\
C_4 \ P_{24} \ ,
\quad
\tr_{(3)}(D_3\ \hat{R}_{23}^{\mp 1}\ \hat{R}_{03}^{\pm 1})\ =\
D_0\ P_{02}\  .
\ee
Consider the left one of equalities in (\ref{tr-04a}) with upper/lower signs.
Multiply both its sides  by $\hat{\Psi}_{23}$/$\hat{\Psi}_{43}$ and take
trace in the space with label 2/4. Then, apply definition (\ref{skew}) and 
use the relation $Tr_{(2)} (U_2 P_{12} W_2) = W_1 U_1$ which  holds for any 
$U,W\in {\rm End}(V)$ and follows from the properties of the trace and the 
permutation. The resulting equality is just the left/right formula 
in (\ref{psi-C}).

Derivation of relations (\ref{psi-D}) 
from the right equality in (\ref{tr-04a}) proceeds similarly.
\hfill$\Box$

\vspace{3mm}

\nin
{\bf Corollary.} Evaluating  traces of relations 
(\ref{psi-C})/(\ref{psi-D}) in spaces with labels 2/1
one finds
\be
\lb{trcd}
Tr_{(2)} C_2 \hat{R}_{21}^{-1} =
Tr_{(2)} D_2 \hat{R}_{12}^{-1} = C_1 D_1 = D_1 C_1 \; .
\ee

\nin
{\bf Theorem.~}
{\it
Let $\hat{R}$ be a skew invertible  BMW type R-matrix. Then the rank of 
the operator $\hat{K}\in \mbox{End}(V^{\otimes 2})$  (see eq.(\ref{imply})) 
equals 1.}
\vspace{2mm}

\nin
{\bf Proof.~}
Consider an R-matrix version of the left equation in (\ref{bmw4})
$$
\hat{K}_{23} \hat{R}_{12} \hat{K}_{23}\ =\ \nu^{-1} \hat{K}_{23} \, .
$$
Multiplying by $\hat{\Psi}_{01}$ and taking trace in space 1 we obtain
\be
\lb{kpk}
\hat{K}_{23} P_{02} \hat{K}_{23}\ =\ \nu^{-1} D_0 \hat{K}_{23} .
\ee
Evaluating traces in spaces 2 and 3 in the left hand side of
relation (\ref{kpk}) one finds
\be
\lb{kpk1}
\tr_{(23)}\left(\hat{K}_{23}P_{02}\hat{K}_{23}\right) =
\mu \ \tr_{(23)}(\hat{K}_{23}P_{02}) =
\mu \ \tr_{(23)}(P_{02}\hat{K}_{03}) =
\mu \ \tr_{(3)}\hat{K}_{03}\ ,
\ee
where (\ref{knorm}) and the properties of the permutation 
were used. On the other hand, (\ref{knorm}) implies
\be
\lb{trK}
\tr_{(12)} \hat{K}_{12} = \mu \ {\rm rank} (\hat{K}) \, ,
\ee
and so, applying  $\tr_{(23)}$ to the right hand side of relation
(\ref{kpk}), one gets
\be
\lb{kpk2}
\nu^{-1} D_0\ \tr_{(23)}\hat{K}_{23}\ =\
{\mu \over \nu}\ {\rm rank} (\hat{K})\ D_0 .
\ee
Equating the results of calculations in (\ref{kpk1}) and (\ref{kpk2})
we obtain
\be
\lb{prep-D-G}
\tr_{(2)} \hat{K}_{12} \ =\ \nu^{-1} {\rm rank} (\hat{K})\ D_1\, .
\ee
In the same way the equality
\be
\lb{prep-C-G}
\tr_{(1)} \hat{K}_{12} \ =\ \nu^{-1} {\rm rank} (\hat{K})\ C_2\ .
\ee
follows from relation (\ref{bmw4b}) with the upper choice of signs.

\medskip
Consider now an R-matrix version of
the right formula in (\ref{bmw4})
$$
\hat{K}_{23} \hat{R}^{-1}_{12} \hat{K}_{23}\ =\ \nu \hat{K}_{23}\, .
$$
Taking traces of this equality in spaces 2 and 3 
and using eq.(\ref{prep-D-G})  one obtains
\be
\lb{Dtr-R-inv}
\tr_{(2)} D_2 \hat{R}_{12}^{-1}\ =\ \nu^2\ I_1\ ,
\ee
which in view of (\ref{psi-D}) is equivalent to
\be
\lb{CD}
CD = DC = \nu^2  I\ .
\ee
Thus, in the conditions of the theorem, the matrices $C$ and $D$
are invertible and we can write relation (\ref{prep-C-G})
in a form
\be
\lb{prepCG}
\tr_{(1)} D_2 \hat{K}_{12} \ =\ \nu \, {\rm rank} (\hat{K})\, I_2\ .
\ee 

On the other hand, applying $\tr_{(23)}$ to an R-matrix version
of eq.(\ref{bmw5}), that is 
$
\hat{K}_{23} \hat{K}_{12} \hat{K}_{23} = \hat{K}_{23} ,
$
and taking into account relations (\ref{trK}) and (\ref{prep-D-G})
we obtain 
\be
\lb{Dtr-G}
\tr_{(2)} D_2 \hat{K}_{12}\ =\ \nu\ I_1\ .
\ee
Finally, evaluating $\tr_{(2)}$ of the equality (\ref{prepCG}) and 
$\tr_{(1)}$ of the equality (\ref{Dtr-G}) and comparing the results we 
conclude that ${\rm rank}\ \hat{K} \ =1$.
\hfill$\Box$

\medskip
\nin {\bf Remark.~} Although in this paper the matrices $C$ and $D$ are 
auxiliary, they play a conceptual role in the theory of quantum groups and 
are used, in particular, for the definition of quantum traces (more details 
on that can be found in \cite{D2,R,O,IOP}). While proving the theorem we 
have derived a number of formulas --- (\ref{Dtr-G}), (\ref{CD}), and 
(\ref{prep-D-G}), (\ref{prep-C-G}) (where one has to put ${\rm rank} 
{\hat K} = 1$) --- which are characteristic for matrices $C$ and $D$ 
corresponding to BMW type R-matrices. One more relation can be added
\be
\lb{trCD}
\tr\, D\ =\ \tr\, C \ =\ \nu \, \mu \ .
\ee
It follows by evaluation of traces of relations (\ref{prep-D-G}) and 
(\ref{prep-C-G}).

\bigskip
>From now on we shall fix some basis  $\{v^i\}_{i=1}^N$ in space $V$ 
($N:=\dim V$). Let
\be
\lb{rank1}
\hat{K}_{ij}^{kl} = \bar{g}_{ij}\ g^{kl}\, .
\ee
be the  matrix of the rank one  operator $\hat{K}$ in this  basis. Define 
operators $X, Y\in {\rm End}(V)$ whose matrices in the chosen basis are 
\be
\lb{XY}
X_i^j := \sum_k g^{ik} \bar{g}_{kj} \, ,
\qquad Y_i^j := \sum_k g^{kj} \bar{g}_{ik} \; .
\ee

\nin {\bf Proposition~2.~} {\it Let $\hat{R}$ be a skew invertible BMW type
R-matrix. Bivectors $g^{kl}$ and $g_{ij}$ (\ref{rank1}) define
nondegenerate bilinear pairings $g: V\otimes V\rightarrow {\Bbb C}$ and 
$\bar{g}:  V^*\otimes V^*\rightarrow {\Bbb C}$ 
\be
\lb{spariv}
g(x,y):=\sum_{i,j=1}^N x_i\, y_j\, g^{ij}\, , \quad
\bar{g}(z,t):=\sum_{i,j=1}^N z^i\, t^j\, \bar{g}_{ij}\, , \quad \forall\;
x,y\in V , \; z,t\in V^*\, ,
\ee
where $x_i$, $y_j$, and $z^i$, $t^j$ stand for coordinates
of vectors $x$, $y$, and $z$, $t$ in the basis $\{v^i\}$ and the dual basis
$\{v^*_i\}$, respectively.
  
Operators $X$ and $Y$ are inverse to each other; the coefficients of the 
characteristic polynomial of $X$: 
$\det(xI-X) = \sum_{k=0}^N(-1)^k C_k\, x^{N-k}\,$ satisfy reciprocity relations
\be
\lb{vozvr}
C_k\, =\, \epsilon\, C_{N-k} \quad (\forall\; 0\leq k\leq N)\, ,\qquad 
\epsilon=\pm 1\, .
\ee
}
\nin
{\bf Proof.~} Consider an R-matrix version of relation (\ref{bmw5}), 
$\hat{K}_{12} \hat{K}_{23} \hat{K}_{12}\ =\ \hat{K}_{12}$.
By a substitution of eq.(\ref{rank1}) and by evaluation of traces in spaces
1 and 2 (note that
$\sum_{i,j}(g^{ij}  \bar{g}_{ij}) = \tr_{(12)}{\hat K}_{12} =  \mu \neq 0$) 
the above equality acquires a form\footnote{
It is this equation which was used in a classification of quantum 
groups in dimension 2 \cite{EOW}.}
\be
\lb{XinvY}
XY =  I\, ,
\ee
wherefrom it also follows that pairings (\ref{spariv}) are nondegenerate.

The definition of matrices $X$ and $Y$ implies that $\tr (X^k) = \tr (Y^k), 
\; \forall k=1,2,\dots $, and, hence, matrices $X$ and $Y=X^{-1}$ obey the 
same characteristic polynomial. Taking into account the identities
$~C_N(X)\, C_k(X^{-1}) = C_{N-k}(X)\, ,$ we then conclude
\be
\lb{zakl}
C_N(X)\, C_k(X)\, =\, C_{N-k}(X)\, \quad \forall\; 1\leq k\leq N\, .
\ee 
For $k=N$ this gives $C_N(X)=\epsilon=\pm 1$; substituting the expression
for $C_N(X)$ back to (\ref{zakl}) one obtains (\ref{vozvr}).
\hfill$\Box$
\vspace{3mm}

Following \cite{FRT} for any  R-matrix $\hat{R}$ we define an associative 
unital bialgebra ${\cal F}(\hat{R})$ generated by components of the matrix
$T:= ||T_i^j||_{i,j=1}^{\phantom{i,j=}N}~$ subject to  relations
\be
\lb{RTT}
\hat{R}_{12} \ T_1 \ T_2 = T_1 \ T_2 \ \hat{R}_{12} \; .
\ee
The coproduct and the counit in the bialgebra are defined by
$$
\triangle (T_i^j) = \sum_{k=1}^N T_i^k \otimes T_k^j \; , \;\;\;
\varepsilon(T_i^j) = \delta_i^j \; .
$$

For the skew invertible BMW type R-matrix $\hat{R}$ the eqs.(\ref{imply}), 
(\ref{RTT}) and the rank one property of the matrix $\hat{K}$
together imply
\be
\lb{kap}
\hat{K}_{12}T_1 T_2\, =\,
\mu^{-1}\hat{K}_{12}T_1 T_2\hat{K}_{12}\, =\,
\tau \hat{K}_{12}\, 
\ee
for some $\tau\in {\cal F}(\hat{R})$.

The following  proposition demonstrates the role of the matrix $X$ for the 
algebra ${\cal F}(\hat{R})$.
\vspace{2mm}

\nin {\bf Proposition~3.~} {\it
Under the assumptions of the theorem
the element $\tau$ is group-like, i.e., $\triangle (\tau) = \tau\otimes\tau$,
$\varepsilon(\tau) = 1$.  It also satisfies relations
\be
\lb{kap1}
\tau\ T_i^j = (XTX^{-1})_i^j\ \tau \; .
\ee
}
\nin
{\bf Proof.~} The group-like properties of the element $\tau$ are directly 
checked by application of the coproduct and the counit operations to 
relation (\ref{kap}). Relation (\ref{kap1}) is justified by a calculation
$$
\tau\ \hat{K}_{12} T_3 = \hat{K}_{12} T_1 T_2 T_3 =
\hat{K}_{12} \hat{R}_{23} \hat{R}_{12} T_1 T_2 T_3
\hat{R}^{-1}_{12} \hat{R}^{-1}_{23} =
$$
$$
\hat{K}_{12} T_1 \hat{K}_{23} T_2 T_3 \hat{R}^{-1}_{12} \hat{R}^{-1}_{23}
= \hat{K}_{12} T_1 \hat{K}_{23} \hat{R}^{-1}_{12} \hat{R}^{-1}_{23}\ \tau =
$$
$$
= \hat{K}_{12} T_1 \hat{K}_{23}
\hat{K}_{12}\ \tau = \hat{K}_{12}
\hat{K}_{23} \hat{K}_{12} (XTX^{-1})_3 \ \tau
= \hat{K}_{12} (XTX^{-1})_3\ \tau \, .
$$
Here we have used relations (\ref{bmw3}), (\ref{bmw3a}), (\ref{bmw5}), 
as well as formula
$$
T_1 {\hat K}_{23} {\hat K}_{12}\, =\,
{\hat K}_{23} {\hat K}_{12} (XTX^{-1})_3\, ,
$$
which is checked by a substitution of expressions 
(\ref{rank1}), (\ref{XY}) for  $\hat{K}$ and $X$
with a subsequent use of (\ref{XinvY}). 
\hfill$\Box$
\vspace{3mm}

Given a quantum group $U_q({\goth g})$ (recall: $\goth g$ is an orthogonal or 
symplectic Lie algebra), its dual  Hopf algebra $G_q$ can be constructed as 
a quotient of the bialgebra ${\cal F}(\hat{R})$ (here the R-matrix $\hat R$ 
is defined by a canonical element of $U_q({\goth g})$) by an ideal $\tau =1$. 
By duality, the left $U_q({\goth g})$-module $V$ admits the right coaction 
of $G_q$
\be
\lb{kodei}
\delta(v^i) = \sum_{j=1}^N v^j\otimes T_j^i\, .
\ee
As one can see from eq.(\ref{kap}) it is the condition $\tau = 1$ that makes 
the pairings $g$ and $\bar g$ (\ref{spariv}) invariant with respect to the 
coaction (\ref{kodei}).

One can start with an algebra ${\cal F}(\hat{R})$ defined by some 
skew invertible BMW type R-matrix $\hat R$. Then the factorization of 
${\cal F}(\hat{R})$ by the relation $\tau =1$ would imply linear dependencies 
among generators $T^i_j$ (c.f., eq.(\ref{kap1})) unless the matrix $X$ 
(\ref{XY}) is scalar. As it is seen from an example below this is not always 
the case.

\medskip
The standard $so_N$ and $sp_N$ series of BMW type R-matrices 
(see \cite{FRT})\footnote{To realize a representation of the
Birman-Murakami-Wenzl algebra the R-matrices given in \cite{FRT}
are to be multiplied (in our case, from the left) by the permutation operator.
} are
\be
\lb{rst}
\hat{R} =\  \sum_{ i,j=1}^N q^{(\delta_{ij}-\delta_{ij'})} \ e_{ij}\otimes 
e_{ji}\ +\ \lambda \sum_{1\leq j<i}^N e_{jj}\otimes e_{ii}
\ -\ \lambda \sum_{1 \leq j<i}^N q^{(\rho_i-\rho_j)}\epsilon_i\epsilon_j\
 e_{i'j}\otimes e_{ij'}\, .
\ee
Here the following notation is used: $i' := N+1-i$;
$||e_{ij}||_k^l:=\delta_{ik}\delta_j^l$ are matrix units;
$\epsilon_i = 1\; (\forall\, i)$ for the ${  s  o}_N$ case
and $\epsilon_i=1=-\epsilon_{i'}\; \forall\; i\leq n$ for the 
${  s  p}_{2n}$ case; the  numbers  
$(\rho_1, \rho_2, \dots ,\rho_N)$ are chosen as 
$(n-1/2, n-3/2, \dots ,1/2,0,-1/2, \dots ,$ $ -n+1/2)$,
$(n-1, n-2, \dots ,1,0,0,-1, \dots ,-n+1)$ and
$(n,$ $n-1,\dots ,1,-1,\dots ,-n)$ in cases   $so_{2n+1}$, $so_{2n}$, or
$sp_{2n}$ correspondingly. In calculations below we will use the relation 
$\rho_i=-\rho_{i'}$ rather then the explicit expressions for
$\rho_i$. The R-matrices (\ref{rst}) generate representations of the 
Birman-Murakami-Wenzl algebras ${  W}_k(q,\nu)$ with specific values of 
their parameter $\nu$, namely, $\nu = q^{1-N}$ for the $so_N$ case and
$\nu = -q^{-1-2N}$ for the $sp_N$ case.

For the R-matrices (\ref{rst}) one calculates
$\bar{g}_{ij} = \delta_{ij'} \epsilon_i q^{-\rho_i}$,
$g^{ij} = \delta^{ij'} \epsilon_{i'} q^{-\rho_i}$, wherefrom it follows 
that $X=I$. To construct R-matrices whose corresponding matrices $X$ 
are not scalars we apply the twist procedure suggested in \cite{R2}
(see also \cite{IOP,Isa}). Remind briefly that given a pair of R-matrices
$\hat R$ and $\hat F$ one can produce a new R-matrix 
$\hat{R}_{F}:= (P{\hat F})\, \hat{R} \, ({\hat F}^{-1}P)$,
called the twisted $\hat R$, provided that additional relations on $\hat R$ 
and $\hat F$ are satisfied
\be
\lb{twist}
\hat{R}_{12} \ \hat{F}_{23} \ \hat{F}_{12} =
\hat{F}_{23} \ \hat{F}_{12} \ \hat{R}_{23} \; , \;\;\;
\hat{F}_{12} \ \hat{F}_{23} \ \hat{R}_{12} =
\hat{R}_{23} \ \hat{F}_{12} \ \hat{F}_{23} \; .
\ee
By construction the twist procedure preserves not only the Yang-Baxter
equation (\ref{YBE}) but all the additional relations 
(\ref{bmw1})--(\ref{bmw4b}) which characterize BMW type R-matrices.

Now we twist R-matrices (\ref{rst}). As a trial twisting R-matrix 
we use $\hat F$ such that $P{\hat F}= \sum_{i,j} d_{ij}\ e_{ii}\otimes 
e_{jj}\ $ where $d_{ij}\in {\Bbb C}\backslash \{0\}$.
An easy check gives the conditions
$$
d_{ij}\ d_{i'j}= u_j\ ,\quad d_{ij}\ d_{ij'}= w_i\ ,\quad \forall\; i,j \, ,
$$
under which relations (\ref{twist}) are satisfied.
The latter in turn are consistent if
$$
u_i u_{i'} = w_i w_{i'} = const \, , \quad \forall \; i\, .
$$
The twisting procedure results in a usual family of multiparametric
R-matrices (some of parameters here are inessential and can be removed by a
linear change of basis in the space $V$)
\be
\lb{multi-R}
\hat{R}_F =\  \sum_{ i,j=1}^N q^{(\delta_{ij}-\delta_{ij'})}
{d_{ij}\over d_{ji}}\ e_{ij}\otimes e_{ji}
\ +\ \lambda \sum_{1\leq j<i}^N e_{jj}\otimes e_{ii}
\ -\ \lambda \sum_{1\leq j<i}^N q^{(\rho_i-\rho_j)}\epsilon_i\epsilon_j\
{d_{i'i}\over d_{jj'}}\ e_{i'j}\otimes e_{ij'} \ .
\ee
For these twisted R-matrices we have
$\bar{g}_{ij} = \delta_{ij'} \epsilon_i q^{-\rho_i} d_{ii'}$,
$g^{ij} = \delta^{ij'} \epsilon_{i'} q^{-\rho_i} d^{-1}_{ii'}$
which gives $X_i^j = \delta_i^j d_{i'i} d^{-1}_{ii'}$.
Thus, element $\tau$ (\ref{kap}) is not necessarily
central in the algebra ${\cal F}(\hat{R})$.
\bigskip

\nin
{\bf Acknowledgements.}~ 
The work of A. P. Isaev and P. N. Pyatov was supported in part by 
the grant No. 03-01-0078 of the Russian Foundation for Basic Research.
The work of A. P. Isaev was also supported by the INTAS grant 
No. 03-51-3350.

\end{document}